\documentclass[12pt,a4paper]{article}
\usepackage{amssymb}

\newtheorem{Thm}{Theorem}

\newtheorem{Cor}[Thm]{Corollary}

\def\RR+{{\mathbf R}^*}

\def\Q_p{{\mathbf Q}_p}

\title{Proper actions of lamplighter groups associated with free groups}

\author{Yves de Cornulier, Yves Stalder and  Alain Valette\footnote{This research was done at Centre Bernoulli (EPF Lausanne), in the framework of the semester ``Limits of graphs in group theory and computer science''.}}

\begin{document}

\baselineskip=16pt

\maketitle

\begin{abstract}
Given a finite group $H$ and a free group $\mathbf F_n$, we prove that the
wreath product $H\wr\mathbf F_n$ admits a metrically proper, isometric action on a Hilbert space.
\end{abstract}

\section{Introduction}

Following \cite{HagPau}, a {\it space with walls} is a pair $(X,{\cal W})$
where $X$ is a set and ${\cal W}$ is a family of partitions of $X$ into two
classes, called {\it walls}, such that for any two distinct points $x,y\in X$,
the number $w(x,y)$ of walls separating $x$ from $y$, called the wall distance between $x$ and $y$, is finite.


Let us define the class ${\cal PW}$ as the class of countable groups $G$
admitting a left-invariant structure of space with walls such that the corresponding wall distance is proper, i.e. bounded subsets are finite. The class ${\cal PW}$ contains for instance $\mathbf{Z}^n$, free groups, surface groups (see 1.2.7 in \cite{CCJJV} for this fact). More generally, it contains all groups acting combinatorially properly on some finite product of trees.

 It is known (see \cite{CCJJV}, Corollary 7.4.2) that, if $G$ belongs to the class ${\cal PW}$, then $G$ has the {\it Haagerup property} (or is
{\it a-T-menable}), i.e. $G$ admits an isometric action on a
Hilbert space $\mathcal{H}$ that is metrically proper, that is, $$\lim_{g\rightarrow\infty} \|x_0 - g\cdot x_0\|=+\infty$$ \nobreak
for some/every $x_0$ in $\mathcal{H}$.

\goodbreak


Recall that the (standard, restricted) {\it wreath product} $H\wr G$ of two groups $H$ and $G$ is the semidirect product $H^{(G)}\rtimes G$, where $G$ acts by shifting the direct sum $H^{(G)}$ of copies of $H$. Up to now, the problem of stability of the Haagerup Property under wreath products was open. Indeed, the Haagerup Property is closed under direct sums, but not under general semidirect products. However it is known to be closed under extensions with {\it amenable} quotients \cite[Example~6.1.6]{CCJJV}, so that in particular $H\wr G$ is Haagerup whenever $H$ is Haagerup and $G$ is amenable.

These were however the only known examples of wreath products with the Haagerup Property. This note presents the first examples of a-T-menable wreath products $H\wr G$ with $H\neq 1$ and $G$ non-amenable, including the first natural such example, namely the ``lamplighter group" $(\mathbf{Z}/2\mathbf{Z})\wr \mathbf{F}_2$ over the free group $\mathbf{F}_2$ of rank two.

\begin{Thm} Let $H$ be a finite group. If $G$ is a group in ${\cal PW}$, then so is $H\wr G$.
In particular, $H\wr G$ has the Haagerup property.\label{th}
\end{Thm}

This latter statement will be generalized in a forthcoming paper, where we prove that the Haagerup Property is closed under taking wreath products. The proof of the general case relies on the same basic construction, but requires more technicalities.


Theorem \ref{th} is proved in Section \ref{proof}, while Section \ref{op} presents some consequences on the relation between the Haagerup Property and weak amenability.

{\bf Ackowledgements:} Thanks are due to I. Chatterji, C. Drutu and S. Popa for useful conversations and correspondence, and to N. Ozawa for suggesting Corollary 3. The
second-named author is especially grateful to Fr\'ed\'eric Haglund for very fruitful
discussions and hints about the construction of families of walls for wreath products.

\section{Proof of Theorem \ref{th}}\label{proof}

We first fix some notation. Write $\Lambda=H^{(G)}$ and $\Gamma= H\wr G= \Lambda\rtimes G$. Elements of $\Gamma$ are denoted
$$\gamma=\lambda g\;(\lambda\in\Lambda, g\in G).$$
The {\it support} of $\lambda$ is
$$\textnormal{supp}(\lambda)=\{g\in G:\lambda(g)\neq 1_H\}.$$

If $(X,{\cal W})$ is a space with walls, a {\it half-space} in $X$ is one of the two classes of some wall in ${\cal W}$. Suppose that $G$ belongs to the class ${\cal PW}$, and let us denote by ${\cal H}$ the family\footnote{Although the family is not assumed injective, we will identify, to avoid subscripts in the notation, elements of the index set ${\cal H}$ with the corresponding subsets of $G$.} of half-spaces in $G$. For $A\in {\cal H}$, we denote by $A^c$ the other half-space in the same wall, i.e. the complement of $A$ in $G$. For $A\in {\cal H}$ and $\mu:A^c\rightarrow H$ a function with finite support, we set
$$E(A,\mu)=:\{\gamma=\lambda g\in\Gamma: g\in A,\;\lambda|_{A^c}=\mu\}.$$

We define a family of walls in $\Gamma$ as partitions $\{E(A,\mu),E(A,\mu)^c\}$, for $A\in{\cal H}$ and $\mu:A^c\rightarrow H$ finitely supported. We check in three steps that this equips $\Gamma$ with a structure of space with walls on which $\Gamma$ acts properly.
\medskip

\underline{1st step}: $\Gamma$ is a space with walls.

Let $\gamma_1=\lambda_1 g_1$ and $\gamma_2=\lambda_2 g_2$ be two elements of $\Gamma$. Let us show that there are finitely many $E(A,\mu)$'s such that $\gamma_1\in E(A,\mu)$ and $\gamma_2\notin E(A,\mu)$.

Indeed $\gamma_1\in E(A,\mu)$ means $g_1\in A$ and $\lambda_1|_{A^c}=\mu$ (so that $\mu$ is determined once $A$ is given). And $\gamma_2\notin E(A,\mu)$ means that either $g_2\notin A$, or $\lambda_2|_{A^c}\neq\mu$; since $\mu=\lambda_1|_{A^c}$, this can be re-written:
$$A^c\cap (\{g_2\}\cup \textnormal{supp}(\lambda_1^{-1}\lambda_2))\neq\emptyset.$$
So $A$ must separate $g_1$ from the finite set $\{g_2\}\cup \textnormal{supp}(\lambda_1^{-1}\lambda_2)$. Since $G$ is a space with walls, this singles out finitely many possibilities for $a$.
\medskip

\underline{2nd step}: $\Gamma$ preserves the above wall structure.

This follows immediately from the formulae:
$$gE(A,\mu)=E(gA,g\mu)\;\;(g\in G);$$
$$\lambda E(A,\mu)=E(a,\lambda|_{A^c}\mu)\;\;(\lambda\in\Lambda).$$
\medskip

\underline{3rd step}: $\Gamma$ acts metrically properly on its wall structure.

Let $w_{\Gamma}(\gamma)$ be the number of walls separating the unit $1_{\Gamma}$ from $\gamma\in\Gamma$. We must prove that, for every $N\in\mathbf{N}$, there are finitely many $\gamma$'s with $w_{\Gamma}(\gamma)\leq N$.

Define analogously $w_G(g)$ as the number of walls separating $1_G$ from $g$ in $G$, and set $B_G(N)=\{g\in G:w_G(g)\leq N\}$; by our assumption $B_G(N)$ is a finite set.

{\bf Claim:} If $\gamma=\lambda g$ satisfies $w_{\Gamma}(\gamma)\leq N$, then $\{g\}\cup \textnormal{supp}(\lambda)\subset B_G(N)$.

Theorem 1 then follows from the claim together with the fact that $H$ is a finite group.

{\bf Proof of the claim:} Contraposing, suppose that there exists $g'\in \{g\}\cup \textnormal{supp}(\lambda)$ with $w_G(g')>N$. So we find $N+1$ distinct half-spaces $A_0,...,A_N$ in $\mathcal{H}$ with $1_G\in A_i$ and $g'\notin A_i\;(i=0,...,N)$. Then the $E(A_i,1_{\Lambda})$'s are $N+1$ distinct half-spaces in $\Gamma$ separating $1_{\Gamma}$ from $\gamma=\lambda g$, so $w_{\Gamma}(\gamma)>N$.
\hfill$\square$

\section{Weak amenability \`a la Cowling-Haagerup}\label{op}

Theorem 1 has interesting consequences in view of a recent result of Ozawa and
Popa \cite{OzPo}. Recall from \cite{CowHaa} that a countable group $G$ is {\it weakly amenable} if there exists a constant $L>0$ and a sequence $(f_n)_{n>0}$ of functions with finite support on $G$, converging pointwise to 1, and such that
$\|f_n\|_{cb}\leq L$ for $n>0$, where $\|f\|_{cb}$ is the
Herz-Schur multiplier norm of the function $f$. The best (i.e. lowest) possible $L$ for which there exists such a sequence is the {\it Cowling-Haagerup constant} of $G$, denoted by $\Lambda(G)$. We set $\Lambda(G)=\infty$ if $G$ is not weakly amenable. Groups which are weakly amenable with constant 1 are also said to satisfy
the {\it complete metric approximation property}.

In \cite[Corollary 2.11]{OzPo}, it is proved that, if $H$ is non-trivial and
$G$ is non-amenable, then $H\wr G$ does not have the complete metric
approximation approximation property. Combining with Theorem 1, we get:

\begin{Cor}\label{cor2}
 For $H$ a non-trivial finite group, $H\wr\mathbf{F}_2$ is an a-T-menable group without the complete metric approximation property.
\hfill$\square$
\end{Cor}

This disproves a conjecture of Cowling (see page 7 in \cite{CCJJV}), stating
that the class of a-T-menable groups coincides with the class of groups with
the complete metric approximation property. Whether every such group is
a-T-menable, is still an open question.

It was pointed out to us by N. Ozawa that from Corollary \ref{cor2} one can deduce the following:

\begin{Cor} Let $H$ be a non-trivial finite group. The iterated wreath product $(H\wr\mathbf{F_2})\wr\mathbf{Z}$ is a-T-menable but not weakly amenable.
\end{Cor}

\noindent {\it Proof:} It was already observed above that $G=:(H\wr\mathbf{F_2})\wr\mathbf{Z}$ is a-T-menable; on the other hand, for $N\geq 1$ consider a subgroup $K_N$ of $G$ which is the direct sum of $N$ copies of $H\wr\mathbf{F}_2$. Then by Proposition 1.3.(a) and Corollary 1.5 in \cite{CowHaa}, we have $\Lambda(G)\geq \Lambda(K_N)=\Lambda(H\wr\mathbf{F_2})^N$. Since $N$ is arbitrary and $\Lambda(H\wr\mathbf{F_2})>1$ by the Ozawa-Popa result \cite{OzPo}, we get $\Lambda(G)=\infty$.\hfill$\Box$

\medskip

In view of the cubulation of spaces with walls, carried out independently in \cite{ChNi} and \cite{Nica}, we get from Theorem 1:

\begin{Cor} Let $H$ be a non-trivial finite group. The wreath product $H\wr\mathbf{F}_2$ admits a metrically proper, isometric action on a $CAT(0)$ cube complex, but does not have the complete metric approximation property.
\hfill$\square$
\end{Cor}

In contrast, it was recently proved by Guentner and Higson \cite{GuHi} that a group acting metrically properly, isometrically on a {\it finite-dimensional} $CAT(0)$ cube complex, has the complete metric approximation property.

\vspace{5mm}

\bibliographystyle{alpha}
\newcommand{\etalchar}[1]{$^{#1}$}

\paragraph{Authors addresses:}
\begin{description}
\item Y.C. Institut de Recherche Math\'ematique de Rennes, Universit\'e de Rennes~1,
Campus de Beaulieu, 35042 Rennes Cedex, France \\
yves.decornulier@univ-rennes1.fr

\item Y.S. Laboratoire de Math\'ematiques, Universit\'e Blaise Pascal, Campus
universitaire des C\'ezeaux, 63177 Aubi\`ere Cedex, France \\
yves.stalder@math.univ-bpclermont.fr

\item A.V. Institut de Math\'ematiques, Universit\'e de Neuch\^atel, Rue \'Emile Argand
11, CP 158, 2009 Neuch\^atel, Switzerland \\
alain.valette@unine.ch

\end{description}
\end{document}